\documentclass{article}             
\usepackage{amssymb}

\def\bR{{\mathbb{R}}}

\def\fg{{\mathfrak{g}}}
\def\fh{{\mathfrak{h}}}
\def\fm{{\mathfrak{m}}}

\def\bA{{\mathbb{A}}}
\def\bB{{\mathbb{B}}}
\def\bC{{\mathbb{C}}}

\newtheorem       {theorem}{Theorem}

\newtheorem{prop} [theorem]{Proposition}
\newtheorem{lemma}[theorem]{Lemma}
\newtheorem{cor}  [theorem]{Corollary}

\begin{document}

\begin{center}
{\Large {Geodesic graphs for geodesic orbit\\
 Finsler $(\alpha,\beta)$ metrics on spheres} }
\bigskip

{\large{Teresa Arias-Marco${}^1$ and Zden\v ek Du\v sek${}^2$}}
\bigskip

${}^1$ Department of Mathematics, University of Extremadura,\\
Av. de Elvas s/n, 06006 Badajoz, Spain\\
ariasmarco@unex.es\\

${}^2$ Institute of Technology and Business in \v Cesk\'e Bud\v ejovice\\
Okru\v zn\'\i\ 517/10, 370 01 \v Cesk\'e Bud\v ejovice, Czech Republic\\
zdusek@mail.vstecb.cz
\end {center}
\bigskip

\begin{abstract} 
Invariant geodesic orbit Finsler $(\alpha,\beta)$ metrics $F$
which arise from Riemannian geodesic orbit metrics $\alpha$ on spheres are determined.
The relation of Riemannian geodesic graphs with Finslerian geodesic graphs
proved in a previous work is now illustrated with explicit constructions.
Interesting examples are found such that $(G/H,\alpha)$ is Riemannian geodesic orbit space,
but for the geodesic orbit property of $(G/H,F)$ the isometry group has to be extended.
It is also shown that projective spaces other than ${\mathbb{R}}P^n$
do not admit invariant purely Finsler $(\alpha,\beta)$ metrics.
\end{abstract}
\bigskip

\noindent
{\bf MSClassification:} {53C22, 53C60, 53C30.}

\section{Preliminaries}
A Finsler manifold $(M,F)$ is {\it homogeneous\/}
if it admits a connected Lie group $G$ which acts transitively on $(M,F)$ as a group of isometries.
Fundamental tools for the study of homogeneous manifolds can be found in several monographs,
for example \cite{De} by S. Deng, concerning Finsler geometry.
The setting convenient for the present work was used for example in the recent work \cite{DuNew} by the second author.
To remain self-contained, we shortly recall these concepts here.
Homogeneous manifold $(M,F)$ can be naturally identified
with the {\it homogeneous space\/} $(G/H,F)$, where $H$ is the isotropy group of the origin $p\in M$.
Denote by $\fg$ and $\fh$ be the Lie algebras of $G$ and $H$
and consider the adjoint representation ${\mathrm{Ad}}\colon H\times\fg\rightarrow\fg$ of $H$ on $\fg$.
There exists a {\it reductive decomposition} of the form $\fg=\fm+\fh$,
where $\fm\subset\fg$ is an ${\rm Ad}(H)$-invariant vector subspace. This decomposition is usually not unique.
For a fixed reductive decomposition $\fg=\fm+\fh$, there is the natural identification
of the vector space $\fm$ with the tangent space of $M$ at $p$ induced by the natural projection
$G\rightarrow G/H=M$.

The invariant Finsler metric $F$ on $G/H$ gives the Minkovski norm and its fundamental tensor on $T_pM$.
Using the above identification $\fm\simeq T_pM$, we obtain the ${\mathrm{Ad}}(H)$-invariant
Minkowski norm and the ${\mathrm{Ad}}(H)$-invariant fundamental tensor on $\fm$.
Conversely, if $F_p$ is given on $\fm\simeq T_pM$, we naturally define $F_q$ on any $T_qM$, $q\in M$, using the action of
$\sigma\in G$ such that $\sigma(p)=q$ by the formula
\begin{eqnarray}
\label{f1}
F_q(\sigma_*X) = F_p(X), \qquad X\in T_pM.
\end{eqnarray}
In the following, we work only with the vector space $\fm$ and we denote
the ${\mathrm{Ad}}(H)$-invariant Minkowski norm on $\fm$ just by $F$ and its fundamental tensor on $\fm$ by $g$.

We shall start with an ${\mathrm{Ad}}(H)$-invariant positive scalar product $\alpha$ and an
${\mathrm{Ad}}(H)$-invariant one-form $\beta$ on $\fm$.
Let $\phi\colon (-b_0,b_0)\rightarrow (0,\infty)$ be a smooth function
and let $\|\beta\|_\alpha<b_0$. Define
\begin{eqnarray}
\label{Fab}
F(y) = \sqrt{\alpha(y,y)}\cdot\phi(s), \qquad s=\frac{\beta(y)}{\sqrt{\alpha(y,y)}}.
\end{eqnarray}
We will use a lemma from \cite{CS}, see also \cite{SS}.
\begin{lemma}[\cite{CS}]
\label{ablemma}
The function $F$ defined in formula $(\ref{Fab})$ is a Minkowski norm for any $\alpha$
and $\beta$ with $\|\beta\|_\alpha<b_0$ if and only if the function $\phi(s)$ satisfies
\begin{equation}
\label{abeq}
\phi(s)>0, \qquad \bigl ( \phi(s)-s\phi'(s)\bigr ) + (b^2 - s^2) \phi''(s) > 0,
\end{equation}
where $s$ and $b$ are arbitrary real numbers with $|s|\leq b < b_0$.
\end{lemma}
The ${\mathrm{Ad}}(H)$-invariant
function $F$ defined on $\fm$ by formula (\ref{Fab}) and which satisfies conditions (\ref{abeq})
obviously determines a $G$-invariant Finsler metric on $M=G/H$.
This metric is called homogeneous $(\alpha,\beta)$ metric.
The notion of an $(\alpha,\beta)$ metric naturally covers Randers metrics for $\phi(s)=1+s$,
quadratic metrics for $\phi(s)=(1+s)^2$ and many other important Finsler metrics,
see for example \cite{CS} or \cite{SS} for more information.

Recall that the standard concepts in Finsler geometry can be found
for example in monographs \cite{BCS} by D. Bao, S.-S. Chern and Z. Shen or \cite{De} by S. Deng.
We consider the {\it Chern connection}, which is the unique linear connection
on the pullback vector bundle $\pi^*TM$ over $TM_0$ which is torsion free and almost $g$-compatible.
Here $g$ is the fundamental tensor related with the metric $F$.
It allows to define the derivative along a curve and geodesics in $M$.

There are many works in which homogeneous geodesics were studied in various settings,
in Riemannian geometry, pseudo-Riemannian geometry, affine geometry and Finsler geometry. See for example
\cite{ASS1}, \cite{ASS2}, \cite{DuS2}, \cite{DuCMUC}, \cite{DuAM}, \cite{DKN},
\cite{GoNi}, \cite{KNi}, \cite{La}, \cite{Ni} and the references therein.
A geodesic $\gamma(s)$ through the point $p\in M$ is {\it homogeneous} if it is an orbit
of a one-parameter subgroup of the group $G=I_0(M)$ of isometries.
If $(M,F)$ is identified with the homogeneous space $(G/H,F)$,
there exists a nonzero {\it geodesic vector} $X\in\fg$ such that $\gamma(t)={\rm exp}(tX)(p)$ for all $t\in\bR$.
A homogeneous space $(G/H, F)$ is called a Finsler {\it geodesic orbit space} or just {\it g.o. space},
if every geodesic of $(G/H, F)$ is homogeneous.

A homogeneous manifold $(M,F)$ may admit more presentations as a homogeneous space in the form $(G/H,F)$,
corresponding to various transitive isometry groups.
There are examples such that $(G/H,F)$ is not a g.o. space, but it is a g.o. space in a presentation
$(\widetilde G/\widetilde H,F)$, where $\widetilde G\supset G$.
In a g.o. space $(G/H,F)$, we investigate sets of geodesic vectors which generate all geodesics through a fixed point.
The concept of geodesic graphs comes from the work of J. Szenthe \cite{Sz}, in the affine setting.
Later, it was used by O. Kowalski with coauthors in the Riemannian setting an by the second author in the Finslerian setting,
see the recent works \cite{DuS2}, \cite{DuCMUC}, \cite{DuAM}, \cite{DKN}, for example.
A {\it geodesic graph} is an ${\mathrm{Ad}}(H)$-equivariant map $\xi\colon\fm\rightarrow\fh$
such that $X +\xi(X)$ is a geodesic vector for each $o\neq X\in\fm$.
If the vector $\xi(X)$ is uniquely determined, then the map $\xi$ is ${\mathrm{Ad}}(H)$-equivariant
and we are interested in the algebraic structure of the mapping $\xi$.
If there are more choices for the vector $\xi(X)$, we want to define geodesic graph
in a way that the algebraic structure of the mapping $\xi$ is as simple as possible.
In \cite{KNi}, O. Kowalski and S. Nik\v cevi\'c developed a theory of Riemannian geodesic graphs.
The existence of a linear geodesic graph is equivalent with the natural reductivity of the space $G/H$.
Or, components of the geodesic graph are rational functions $\xi_i=P_i/P$, where $P_i$ and $P$
are homogeneous polynomials and it holds ${\mathrm{deg}}(P_i)={\mathrm{deg}}(P)+1$.

Geodesic graphs on many Riemannian g.o. manifolds were constructed in several works,
see for example the recent survey paper \cite{DuS2} by the author for more details and references
and also for many important classes of Riemannian g.o. spaces.
Further references and a structural approach to Riemannian g.o. manifolds using Lie theory
can be found in the recent papers \cite{GoNi} and \cite{Ni} by C.S. Gordon and Yu.G. Nikonorov.
Recent results about Riemannian g.o. spaces of special types were obtained for example
in \cite{ASS1}, \cite{ASS2} or \cite{CCZ}.
Finsler g.o. spaces attained attention recently.
In \cite{YD}, Z. Yan and S. Deng studied relations of special Finsler g.o. spaces and Riemannian g.o. spaces.
Particular results were obtained for Randers g.o. metrics and for weakly symmetric metrics.
Invariant Randers g.o. metrics on spheres $S^{2n+1}={\mathrm{U}}(n+1)/{\mathrm{U}}(n)$
were also considered in this paper.
In \cite{DuCMUC}, the second author investigated invariant Randers geodesic orbit metrics
on modified H-type groups and he constructed Finslerian geodesic graphs on these geodesic orbit manifolds.
In \cite{DuAM}, the second author investigated special families of weakly symmetric Finsler g.o. metrics
on modified H-type groups which were studied also in \cite{YD}.
The relation of Finslerian geodesic graphs with the Riemannian geodesic graphs
of the underlying Riemannian metric was studied.

In \cite{DuNew}, second author studied general Finsler $(\alpha,\beta)$ metrics which arise from a Riemannian
geodesic orbit metric $\alpha$ and a one-form $\beta$.
In particular, the class of $(\alpha,\beta)$ metrics covers all Randers metrics as a special case.
Geodesic lemma, which is the formula for characterization of geodesic vectors, was derived
in terms of the Riemannian metric $\alpha$ and the one-form $\beta$.
The existence of a suitable reductive decomposition $\fg=\fm+\fh$ was proved in which the relation
of Riemannian geodesic graph with the Finslerian geodesic graph can be easily described,
in a suitable group extension, without the explicit construction of the geodesic graph.
As a consequence, it was proved that all $(\alpha,\beta)$ metrics such that $\alpha$ is a Riemannian
geodesic orbit metric are Finsler geodesic orbit metrics, possibly with respect to a bigger group of isometries.
It was an open question whether this group extension is necessary only for the particular theoretical
construction in the proof there or whether there are examples which require this extension
to obtain the geodesic orbit property for the Finsler metrics.
In the special situation when the Riemannian g.o. metric $\alpha$ is naturally reductive,
an alternative description of Finslerian geodesic graph with respect to the naturally reductive
decomposition was also given, again, without the explicit construction of geodesic graph.
The construction was illustrated with one example which is not naturally reductive.

In the present paper, we are going to illustrate these construction on spheres,
and we observe interesting phenomena.
For the list of Riemannian geodesic orbit metrics on spheres we refer for example for the monograph \cite{BN}
by V. Berestovskii and Yu.~Nikonorov or for the references therein for the detailed results in particular cases.
In each case, we determine invariant Finsler $(\alpha,\beta)$ metrics and analyze the geodesic orbit property.
In each of the series, we work just with the low-dimensional example.
Due to the explicit construction, it can be easily seen that
similar results hold for the same type of spaces in higher dimensions.
For the general case, we argue using Corollary \ref{cor1} below.

Spheres $S^{2n+1}={\mathrm{SU}}(n+1)/{\mathrm{SU}}(n)$ with any geodesic orbit Riemannian metric $\alpha$
with respect to the group ${\mathrm{SU}}(n+1)$ are also geodesic orbit with respect to this group with all
$(\alpha,\beta)$ Finsler metrics which arise from the Riemannian metric $\alpha$.
However, for the easy description of the relation of Riemannian and Finslerian geodesic graph,
we need the extended expression $S^{2n+1}={\mathrm{U}}(n+1)/{\mathrm{U}}(n)$.
But this space with Riemannian g.o. metrics is naturally reductive and it is more convenient to use
the naturally reductive decomposition, hence we illustrate this approach.
For spheres $S^{4n+3}={\mathrm{Sp}}(n+1)/{\mathrm{Sp}}(n)$, there is a $2$-parameter family of Riemannian
geodesic orbit metrics.
There are many invariant $(\alpha,\beta)$ metrics, but none of them is geodesic orbit
with respect to ${\mathrm{Sp}}(n+1)$.
For each invariant $(\alpha,\beta)$ metric, the isometry group can be extended to
${\mathrm{Sp}}(n+1)\cdot {\mathrm{U}}(1)$ and the relation of Riemannian and Finslerian geodesic graph
can be easily described with respect to this group.
The isometry group of all these Riemannian metrics can be further extented to
${\mathrm{Sp}}(n+1)\cdot {\mathrm{Sp}}(1)$ to obtain natural reductivity,
but there are no invariant $(\alpha,\beta)$ metrics with respect to this group.
We also show that the exceptional spheres which admit Riemannian geodesic orbit metrics
do not admit invariant purely Finsler $(\alpha,\beta)$ metrics.
So our study concludes with the classification of geodesic orbit Finsler $(\alpha,\beta)$ metrics on spheres.
We also examine the projective spaces to see that
other projective spaces than those ${\mathbb{R}}P^n$ which arise naturally from the spheres
do not admit invariant Finsler $(\alpha,\beta)$ metrics.

\section{Geodesic graphs}
We recall geodesic lemma, which characterizes geodesic vectors.
For general Finsler metrics it was derived by D. Latifi.
\begin{lemma}[\cite{La}]
\label{golema2}
Let $(G/H,F)$ be a homogeneous Finsler space with a reductive decomposition $\fg=\fm+\fh$.
A nonzero vector $Y\in{\fg}$ is geodesic vector if and only if it holds
\begin{equation}
\label{gl2}
g_{Y_\fm} ( Y_{\mathfrak m}, [Y,U]_{\mathfrak m} ) = 0 \qquad \forall U\in{\mathfrak m},
\end{equation}
where the subscript $\fm$ indicates the projection of vector from $\fg$ to $\fm$.
\end{lemma}
In \cite{DuNew}, the second author expressed this condition for Finsler $(\alpha,\beta)$ metrics
in terms of the Riemannian metric $\alpha$ and the one-form $\beta$.
Recall that the one-form $\beta$ corresponds to an $\alpha$-equivalent vector $V\in\fm$
related with $\beta$ by the condition $\beta(U) = \alpha(V,U)$ for all $U\in \fm$.
We further denote by $\zeta=\zeta(X)$ the following function on $\fm$ characteristic
for the particular $(\alpha,\beta)$ metric. We write $\phi'={\mathrm{d}}\phi/{\mathrm{d}}s$.
\begin{eqnarray}
 \zeta = \zeta(X) = \frac{\alpha\phi'}{\sqrt{\alpha}\phi-\beta\phi'}.
\end{eqnarray}
\begin{lemma}[\cite{DuNew}]
\label{c2}
Let $F=\sqrt{\alpha}\cdot\phi(s)$ be a homogeneous Finsler $(\alpha,\beta)$ metric on $G/H$,
let $\fg=\fm+\fh$ be a reductive
decomposition and $V\in\fm$ be the vector $\alpha$-equivalent with $\beta$.
The vector $X+\xi(X)$, where $X\in\fm$ and $\xi(X)\in\fh$, is geodesic vector if and only if it holds
\begin{eqnarray}
\label{goleman2}
\alpha \Bigl ( X + \zeta(X) \cdot  V,
 [X+\xi(X),U]_{\mathfrak m}\Bigr ) = 0, \qquad \forall U\in\fm.
\end{eqnarray}
\end{lemma}
Let us also recall the algebraic feature which is necessary for the further results.
\begin{prop}[\cite{DuNew}]
\label{prop1}
Let $(G/H,\alpha)$ be a Riemannian geodesic orbit space with a reductive decomposition $\fg=\fh+\fm$
and let $V\in \fm$ be an ${\mathrm{Ad}(H)}$-invariant vector.
Then either there exist an invariant reductive decomposition $\fg=\fh+\fm'$ such that the projection
of the vector $V$ to $\fm'$ is in the center of $\fg$,
or $M=G/H$ can be expressed in an extended form $M=\widetilde G/\widetilde H$ and the above property
holds with respect to a decomposition $\widetilde\fg=\widetilde\fh+\fm'$.
\end{prop}
The crucial step in the proof is the fact that the operator ${\mathrm{ad}}(V)|_{\fm}$ acts
skew-symmetrically on the suitable reductive complement $\fm$ and we can define new formal operator
$W={\mathrm{ad}}(V)|_{\fm}$ and extend the isotropy algebra by this operator, if necessary.
Consequently, the vector $V-W$ generates nontrivial center in $\widetilde\fg$.
The following easy corollary will be crucial later.
\begin{cor}
\label{cor1}
Let $(G/H,\alpha)$ be a Riemannian geodesic orbit space.
The Riemannian metric $\alpha$ can be modified into a Finslerian $(\alpha,\beta)$ metric $F$
if and only if either $\fg$ has nontrivial center or there exists the extended expression of the manifold
in the form $\widetilde G/\widetilde H$, where $\widetilde\fg$ has nontrivial center.
\end{cor}
{\it Proof.}
Projection of any vector from the center of $\widetilde\fg$ into any reductive complement $\fm$
gives an ${\mathrm{Ad}(H)}$-invariant vector $V\in \fm$ and invariant Finsler $(\alpha,\beta)$ metrics,
according to Lemma \ref{ablemma}.
On the other hand, existence of invariant Finsler $(\alpha,\beta)$ metric requires
an ${\mathrm{Ad}(H)}$-invariant vector $V\in \fm$ and according to Proposition \ref{prop1},
there must be nontrivial center in $\widetilde\fg$.
$~\hfill\square$

Proposition \ref{prop1} is also suitable for an easy construction of Finslerian geodesic graph
using the Riemannian one. We shall illustrate these constructions in the next section.
\begin{theorem}[\cite{DuNew}]
\label{t2}
Let $(G/H,\alpha)$ be a Riemannian geodesic orbit space.
Then all invariant Finsler $(\alpha,\beta)$ metrics $F=\sqrt{\alpha}\cdot\phi(s)$ on $G/H$ are geodesic orbit metrics,
possibly with respect to a bigger group of isometries.
The geodesic graph is
\begin{eqnarray}
\xi(X)=\xi_R(X+\zeta(X)\cdot V),
\end{eqnarray}
where $\xi_R$ is the Riemannian geodesic graph
in a suitable group extension and with respect to a suitable reductive decomposition.
\end{theorem}
If the Riemannian geodesic orbit metric $\alpha$ on $M$ is naturally reductive,
the decomposition from Proposition \ref{prop1} is usually not naturally reductive.
An alternative construction in the naturally reductive decomposition is the following.
\begin{prop}[\cite{DuNew}]
\label{pnr}
Let $(G/H,\alpha)$ be naturally reductive Riemannian homogeneous space with the naturally reductive
decomposition $\fg=\fm+\fh$ and with the nontrivial center ${\mathfrak{c}}\subset\fg$.
For any Finsler $(\alpha,\beta)$ metric
$F=\sqrt{\alpha}\cdot\phi(s)$ determined by the Riemannian metric $\alpha$
and an ${\mathrm{Ad}}(H)$-invariant vector $V=C_\fm$, such that $C_\fm+C_\fh\in{\mathfrak{c}}$,
the geodesic graph is
\begin{eqnarray}
\xi(X)=-\zeta(X)\cdot C_\fh.
\end{eqnarray}
\end{prop}

\section{Geodesic orbit $(\alpha,\beta)$ metrics on spheres}
We are going to consider all expressions of spheres as a homogeneous space with Riemannian
geodesic orbit metrics, according to the classification table in \cite{BN}.
We analyze the existence of invariant Finsler $(\alpha,\beta)$ metrics and we construct
geodesic graphs on them to illustrate Theorem \ref{t2} and Proposition \ref{pnr}.
For the infinite series of spheres, we shall describe explicitly the low-dimensional examples.

We start with the standard representation of the sphere
$S^n={\mathrm{SO}}(n+1)/{\mathrm{SO}}(n)$.
It is well known that it admits one-parameter family of invariant Riemannian metrics
and these metrics are naturally reductive. The isotropy representation does not admit
any invariant vectors in $\fm$ and hence no invariant Finsler $(\alpha,\beta)$ merics.
We shall analyze other expressions of some spheres, for smaller isometry groups,
which admit more invariant geodesic orbit metrics and also some invariant vectors in $\fm$
and hence some invariant Finsler $(\alpha,\beta)$ metrics.
Namely, $S^{2n+1}={\mathrm{SU}}(n+1)/{\mathrm{SU}}(n)$ which admits group extension and expression
$S^{2n+1}={\mathrm{U}}(n+1)/{\mathrm{U}}(n)$
and $S^{4n+3}={\mathrm{Sp}}(n+1)/{\mathrm{Sp}}(n)$ which admits two possible group extensions,
see also \cite{BN} for details about Riemannian geodesic orbit metrics.
Finally, we show that exceptional spheres do not admit invariant purely Finsler $(\alpha,\beta)$ metrics
and we classify all Finsler $(\alpha,\beta)$ geodesic orbit metrics on spheres.

\subsection{$S^{2n+1}={\mathrm{SU}}(n+1)/{\mathrm{SU}}(n)$}
\label{31}
For $n=1$, we are reduced to the situation $S^{3}={\mathrm{SU}}(2)$.
In this case, the geodesic orbit Riemannian metrics $\alpha$ are multiples of the bi-invariant metric.
The isometry group can be enlarged and this example can be treated in the framework of the next section.
To illustrate the general behaviour of the present series,
we shall consider $n=2$ and hence $S^5={\mathrm{SU}}(3)/{\mathrm{SU}}(2)$.
On the Lie algebra level, we choose the basis $\{H_1,H_2,H_3\}$ of the Lie algebra $\fh={\mathfrak{su}}(2)$,
given by the matrices
\begin{eqnarray}
\nonumber
H_1= \left ( \begin{array} {ccc}
i & 0 & 0
\\{\medskip}
0 &-i & 0
\\{\medskip}
0 & 0 & 0
\end{array} \right ),
H_2= \left ( \begin{array} {ccc}
0 & i & 0
\\{\medskip}
i & 0 & 0
\\{\medskip}
0 & 0 & 0
\end{array} \right ),
H_3= \left ( \begin{array} {ccc}
0 & 1 & 0
\\{\medskip}
-1 & 0 & 0
\\{\medskip}
0 & 0 & 0
\end{array} \right ).
\end{eqnarray}
For the reductive complement $\fm$ in the decomposition ${\mathrm{su}}(3)= {\mathrm{su}}(2) + \fm$,
we choose the basis $B=\{X_1,X_2,Y_1,Y_2,Z\}$, given by the matrices
\begin{eqnarray}
\nonumber
&
X_1= \left ( \begin{array} {ccc}
0 & 0 & 1
\\{\medskip}
0 & 0 & 0
\\{\medskip}
-1 & 0 & 0
\end{array} \right ), \quad
X_2= \left ( \begin{array} {ccc}
0 & 0 & i
\\{\medskip}
0 & 0 & 0
\\{\medskip}
i & 0 & 0
\end{array} \right ),\cr
&
Y_1= \left ( \begin{array} {ccc}
0 & 0 & 0
\\{\medskip}
0 & 0 & 1
\\{\medskip}
0 & -1 & 0
\end{array} \right ), \quad
Y_2= \left ( \begin{array} {ccc}
0 & 0 & 0
\\{\medskip}
0 & 0 & i
\\{\medskip}
0 & i & 0
\end{array} \right ), \quad
\nonumber
Z= \left ( \begin{array} {ccc}
-\frac{i}{2} & 0 & 0
\\{\medskip}
0 & -\frac{i}{2} & 0
\\{\medskip}
0 & 0 & i
\end{array} \right ).
\end{eqnarray}
By the straightforward calculations, we obtain the Lie bracket relations
on $\fh={\mathrm{span}}\{H_i\}$, which are
\begin {eqnarray}
\nonumber
 [ H_1,H_2 ] = -2 H_3, \quad [ H_1,H_3 ] =  2 H_2, \quad [ H_2,H_3 ] = -2 H_1
\end{eqnarray}
and the Lie bracket relations on $\fm={\mathrm{span}}\{X_i,Y_i,Z\}$, which are
\begin{eqnarray}
\label{z31}
 & & [ X_1,X_2 ] = -2 Z + H_1, \cr
 & & [ X_1,Y_1 ] =  - H_3, \quad [ X_2,Y_1 ] =  - H_2, \cr
 & & [ X_1,Y_2 ] = ~~ H_2, \quad [ X_2,Y_2 ] =  - H_3, \quad [ Y_1,Y_2 ] = -2 Z - H_1, \cr
 & & [ X_1,Z ] = \frac{3}{2} X_2, \quad [ X_2,Z ] =-\frac{3}{2} X_1, \quad
  [ Y_1,Z ] = \frac{3}{2} Y_2, \quad [ Y_2,Z ] =-\frac{3}{2} Y_1.~~~~~~
\end{eqnarray}
The adjoint action of $\fh={\mathfrak{su}}(2)$ on $\fm$ is generated by the operators
\begin{eqnarray}
\label{a31}
{\rm ad}(H_1)|_\fm = A_{12}-A_{34},\quad
{\rm ad}(H_2)|_\fm = A_{14}-A_{23},\quad
{\rm ad}(H_3)|_\fm =-A_{13}-A_{24}.
\end{eqnarray}
From the adjoint action given by formulas (\ref{a31}) we see that
there is a 2-parameter family of invariant scalar products on $\fm$
which determine invariant Riemannian metrics on ${\mathrm{SU}}(3)/{\mathrm{SU}}(2)$.
Up to a multiple, these scalar products are determined by the condition that the basis
$\{X_1,X_2,Y_1,Y_2,\frac{1}{\sqrt{c}}Z\}$, for some $c>0$, is orthonormal.
These Riemannian metrics are known to be weakly symmetric and hence geodesic orbit, see for example
\cite{BN}, \cite{Z1}, \cite{Z2}, \cite{Z3}.
\begin{prop}
Each of the $2$-parameter family of Riemannian geodesic orbit metrics on the sphere
$S^5={\mathrm{SU}}(3)/{\mathrm{SU}}(2)$ can be modified in one direction to obtain
invariant Finsler $(\alpha,\beta)$ metrics. All these metrics are geodesic orbit with respect
to the group ${\mathrm{SU}}(3)$.
\end{prop}
{\it Proof.}
The vector $Z\in\fm$ is ${\mathrm{Ad}}(H)$-invariant and hence its $\alpha$-equivalend one-form $\beta$,
together with the Riemannian metric $\alpha$ and a function $\phi(s)$,
according to Lemma \ref{ablemma}, determine invariant Finsler $(\alpha,\beta)$ metrics.
We are now going to illustrate the geodesic orbit property of these $(\alpha,\beta)$ metrics.
We use equation (\ref{goleman2}) in the form
\begin{eqnarray}
\label{forex}
\alpha \Bigl ( X , [\xi(X),U]_{\mathfrak m}\Bigr ) & = &
- \alpha \Bigl ( X , [X,U]_{\mathfrak m}\Bigr ) 
- \alpha \Bigl ( \zeta(X)\cdot  V, [X,U]_{\mathfrak m}\Bigr ) ,
\end{eqnarray}
which allows us to construct both the Riemannian geodesic graph $\xi_R$ and the Finslerian geodesic graph $\xi$.
We write each vector $X\in\fm$ and each vector $F\in\fh$ in the forms
\begin{eqnarray}
\nonumber
X & = & x_1X_1 + x_2 X_2 + x_3 Y_1 + x_4 Y_2 + z Z, \cr
F & = & \xi_1 H_1 + \dots +\xi_3 H_3,
\end{eqnarray}
and we consider formula (\ref{forex}) to determine components $\xi_i$ depending on $x_j$ and on $z$.
We obtain the system of linear equations and we write down the extended matrix $(\bA,\bB,\bC)$ of this system,
where we separate the two right-hand sides according the two terms on the right-hand side in formula (\ref{forex}).
The second right-hand side, which we denote by $(\bC)$ and which corresponds
to the last term in formula (\ref{forex}) is the Finslerian term.
The extended matrix is
\begin{eqnarray}
\nonumber
(\bA | \bB | \bC) =
\left (
\begin{array}
{ccc|c|c}
{x_2}&{x_4}&{-x_3}&
 x_2z(\frac{3}{2}-2c) & -2 x_2 c \zeta v
\\{\medskip}
{-x_1}&{-x_3}&{-x_4}&
-x_1z(\frac{3}{2}-2c) & 2 x_1 c \zeta v
\\{\medskip}
{-x_4}&{x_2}&{x_1}&
 x_4z(\frac{3}{2}-2c) & -2 x_4 c \zeta v
\\{\medskip}
{x_3}&{-x_1}&{x_2}&
-x_3z(\frac{3}{2}-2c) & 2 x_3 c \zeta v
\end{array}
\right ).
\end{eqnarray}
The rank of this system is equal to 3 and it is solvable by the Cramer's rule.
In the Riemannian case and also in the Finslerian case,
the solution is the unique geodesic graph.
After cancelling out the common factor, its components are $\xi_i=\frac{P_i}{P}$,
where ${\mathrm{deg}}(P_i)-1={\mathrm{deg}}(P)=2$.
It shows that all the Finsler $(\alpha,\beta)$ metrics which arise from the Riemannian
geodesic orbit metric $\alpha$ and the one-form $\beta$ which is $\alpha$-equivalent with the vector $Z\in\fm$
are also geodesic orbit metrics.
$\hfill\square$

We remark that the geodesic graph constructed from the system
$(\bA | \bB | \bC)$ above proves the geodesic orbit property of the space $G/H$ with respect to he group $G$,
but it is not constructed according to Theorem \ref{t2}.

\subsection{$S^{2n+1}={\mathrm{U}}(n+1)/{\mathrm{U}}(n)$}
\label{32}
We continue with the case $n=2$, hence we have $S^5={\mathrm{U}}(3)/{\mathrm{U}}(2)$.
We are going to illustrate Proposition \ref{pnr} for the space from previous section,
in the extended isometry group $\widetilde G={\mathrm{U}}(3)$.
The space with Riemannian metrics $\alpha$ considered above
is known to be naturally reductive, with respect to the isometry group ${\mathrm{U}}(3)$,
see for example \cite{BN}, \cite{Z1}, \cite{Z2}, \cite{Z3}.
We denote further
\begin{eqnarray}
\nonumber
H_0= \left ( \begin{array} {ccc}
i & 0 & 0
\\{\medskip}
0 & i & 0
\\{\medskip}
0 & 0 & 0
\end{array} \right ), \quad
\bar Z= \left ( \begin{array} {ccc}
i(1-2c) & 0 & 0
\\{\medskip}
0 & i(1-2c) & 0
\\{\medskip}
0 & 0 & i
\end{array} \right )
\end{eqnarray}
and consider the basis $\{H_0,H_1,H_2,H_3\}$ of the Lie algebra $\widetilde\fh={\mathfrak{u}}(2)$.
We have ${\mathrm{u}}(3)= {\mathrm{u}}(2) + \fm$, where the subspace $\fm$ is again generated by the basis $B$,
as in the previous section. The adjoint action of $\widetilde\fh$ on $\fm$ is given by the operators
\begin{eqnarray}
\label{a32}
{\rm ad}(H_1)|_\fm = A_{12}-A_{34},\quad
{\rm ad}(H_2)|_\fm = A_{14}-A_{23},\cr
{\rm ad}(H_3)|_\fm =-A_{13}-A_{24},\quad
{\rm ad}(H_0)|_\fm = A_{12}+A_{34}.
\end{eqnarray}
The invariant scalar products on $\fm$ are the same as in previous section and they determine
invariant Riemannian metrics also on ${\mathrm{U}}(3)/{\mathrm{U}}(2)$.
Again, the vector $Z\in\fm$ is ${\mathrm{Ad}}(H)$-invariant and hence,
together with a function $\phi(s)$, according to Lemma \ref{ablemma},
it determines invariant Finsler $(\alpha,\beta)$ metrics.
For $F\in\widetilde\fh$, we write
\begin{eqnarray}
\nonumber
F & = & \xi_1 H_1 + \dots +\xi_3 H_3 + \xi_4 H_0
\end{eqnarray}
and we consider again formula (\ref{forex}).
We obtain the system of linear equations whose extended matrix is
\begin{eqnarray}
\nonumber
(\bA | \bB | \bC) =
\left (
\begin{array}
{cccc|c|c}
{x_2}&{x_4}&{-x_3}&{x_2} &
 x_2z(\frac{3}{2}-2c) & -2 x_2 c \zeta v
\\{\medskip}
{-x_1}&{-x_3}&{-x_4}&{-x_1} &
-x_1z(\frac{3}{2}-2c) & 2 x_1 c \zeta v
\\{\medskip}
{-x_4}&{x_2}&{x_1}&{x_4} &
 x_4z(\frac{3}{2}-2c) & -2 x_4 c \zeta v
\\{\medskip}
{x_3}&{-x_1}&{x_2}&{-x_3} &
-x_3z(\frac{3}{2}-2c) & 2 x_3 c \zeta v
\end{array}
\right ).
\end{eqnarray}
Obviously, there exist a linear geodesic graph $\xi(X)=z(\frac{3}{2}-2c)\cdot H_0$.
We now change into the naturally reductive decomposition $\widetilde\fg=\widetilde\fh+\fm'$, where
$\fm'$ is generated by the basis $B=\{X_i,Y_i,\bar Z\}$ and
$\bar Z=Z+\xi(Z)=Z+(\frac{3}{2}-2c)\cdot H_0$.
The nonzero projections of the new Lie brackets of elements from $\fm'$ into $\fm'$ are
\begin{eqnarray}
\nonumber
 & & [ X_1,X_2 ]_{\fm'} =  -2 \bar Z , \quad [ Y_1,Y_2 ]_{\fm'} = -2 \bar Z , \cr
 & & [ X_1,\bar Z ]_{\fm'} = 2c X_2, \quad [ X_2,\bar Z ]_{\fm'} =-2c X_1, \quad
  [ Y_1,\bar Z ]_{\fm'} = 2c Y_2, \quad [ Y_2,\bar Z ]_{\fm'} =-2c Y_1.
\end{eqnarray}
The adjoint action of $\widetilde\fh={\mathfrak{u}}(2)$ on $\fm'$ is given again by the operators (\ref{a32}).
Formula (\ref{forex}) in the naturally reductive decomposition $\widetilde\fg=\widetilde\fh+\fm'$ is in the form
\begin{eqnarray}
\nonumber
\alpha \Bigl ( X , [\xi(X),U]_{\mathfrak m}\Bigr ) & = &
- \alpha \Bigl ( \zeta(X)\cdot  V, [X,U]_{\mathfrak m}\Bigr ),
\end{eqnarray}
because the first term on the right-hand side vanishes.
The matrix of the corresponding system of equations is now
\begin{eqnarray}
\nonumber
(\bA | \bC) =
\left (
\begin{array}
{cccc|c}
{x_2}&{x_4}&{-x_3}&{x_2} &
 -2 x_2 c \zeta v
\\{\medskip}
{-x_1}&{-x_3}&{-x_4}&{-x_1} &
 2 x_1 c \zeta v
\\{\medskip}
{-x_4}&{x_2}&{x_1}&{x_4} &
 -2 x_4 c \zeta v
\\{\medskip}
{x_3}&{-x_1}&{x_2}&{-x_3} &
 2 x_3 c \zeta v
\end{array}
\right ).
\end{eqnarray}
We see that there is the solution $\xi(X)=-2c\zeta(X)\cdot H_0$.
Because $\bar Z + 2cH_0\in {\mathfrak{c}}(\fg)$,
this solution is the Finslerian geodesic graph,
according to Proposition \ref{pnr}.
\begin{prop}
\label{pr1}
Each of the $2$-parameter family of Riemannian naturally reductive metrics on the sphere
$S^{2n+1}={\mathrm{U}}(n+1)/{\mathrm{U}}(n)$ can be modified in one direction to obtain
invariant Finsler $(\alpha,\beta)$ metrics. All these metrics are geodesic orbit with respect
to the group ${\mathrm{U}}(n+1)$.
\end{prop}
{\it Proof.}
The existence of the $2$-parameter family of Riemannian naturally reductive metrics is well known,
see for example \cite{BN}, \cite{Z1}, \cite{Z2}, \cite{Z3}.
The detailed construction of these metrics and also the statement of the proposition
in the particular case $n=2$ is clear from the example above.
Concerning the general case, we use Corollary \ref{cor1} and Proposition \ref{pnr}.
Let $\fg=\fm+\fh$ be the naturally reductive decomposition 
and ${\mathfrak{c}}\subset\fg={\mathfrak{u}}(n+1)$ be the center.
For any ${\mathrm{Ad}}(H)$-invariant vector $V=C_\fm\in\fm$, there is a unique vector $C_\fh\in\fh$
such that $C_\fm+C_\fh\in{\mathfrak{c}}$.
Geodesic graph for any invariant Finsler $(\alpha,\beta)$ metric $F$ determined by the Riemannian
metric $\alpha$, the one-form $\beta$ which is $\alpha$-equivalent with $V$ and a function $\phi(s)$,
according to Lemma \ref{ablemma}, can be constructed using Proposition \ref{pnr}.
This proves the geodesic orbit property of $(G/H,F)$.
$\hfill\square$

\subsection{$S^{4n+3}={\mathrm{Sp}}(n+1)/{\mathrm{Sp}}(n)$ }
\label{33}
For $n=1$, we have $S^7={\mathrm{Sp}}(2)/{\mathrm{Sp}}(1)$.
On the Lie algebra level, we choose a basis $\{H_1,H_2,H_3\}$ of the Lie algebra $\fh={\mathfrak{sp}}(1)$,
given by the matrices
\begin{eqnarray}
\nonumber
H_1= \left ( \begin{array} {ccc}
i & 0 
\\{\medskip}
0 & 0
\end{array} \right ),
H_2= \left ( \begin{array} {ccc}
j & 0
\\{\medskip}
0 & 0
\end{array} \right ),
H_3= \left ( \begin{array} {ccc}
k & 0
\\{\medskip}
0 & 0
\end{array} \right ).
\end{eqnarray}
For the reductive complement $\fm$ in the decomposition ${\mathrm{sp}}(2)= {\mathrm{sp}}(1) + \fm$,
we choose the basis $B=\{X_1,\dots,X_4,Z_1,\dots,Z_3\}$, given by the matrices
\begin{eqnarray}
\nonumber
&
X_1= \left ( \begin{array} {ccc}
0 & 1
\\{\medskip}
-1 & 0
\end{array} \right ),
X_2= \left ( \begin{array} {ccc}
0 & i
\\{\medskip}
i & 0
\end{array} \right ),
X_3= \left ( \begin{array} {ccc}
0 & j
\\{\medskip}
j & 0
\end{array} \right ),
X_4= \left ( \begin{array} {ccc}
0 & k
\\{\medskip}
k & 0
\end{array} \right ),\cr
&
Z_1= \left ( \begin{array} {ccc}
0 & 0 
\\{\medskip}
0 & i
\end{array} \right ),
Z_2= \left ( \begin{array} {ccc}
0 & 0
\\{\medskip}
0 & j
\end{array} \right ),
Z_3= \left ( \begin{array} {ccc}
0 & 0
\\{\medskip}
0 & k
\end{array} \right ).
\end{eqnarray}
By the straightforward calculations, we obtain the Lie bracket relations on $\fh$, which are
\begin {eqnarray}
\nonumber
 [ H_1,H_2 ] =  2 H_3, \quad [ H_1,H_3 ] = -2 H_2, \quad [ H_2,H_3 ] = 2 H_1.
\end{eqnarray}
The Lie bracket relations on $\fm={\mathrm{span}}\{X_i,Z_j\}$ are
\begin{eqnarray}
\nonumber
 & & [ X_1,X_2 ] = -2 Z_1 + 2 H_1, \cr
 & & [ X_1,X_3 ] = -2 Z_2 + 2 H_2, \quad [ X_2,X_3 ] =  2 Z_3 + 2 H_3, \cr
 & & [ X_1,X_4 ] = -2 Z_3 + 2 H_3, \quad [ X_2,X_4 ] =  -2 Z_2 - 2 H_2, \quad [ X_3,X_4 ] = ~2 Z_1 + 2 H_1,
\end{eqnarray}
\begin{eqnarray}
\nonumber
 & & [ X_1,Z_1 ] =  X_2, \quad [ X_2,Z_1 ] =-X_1, \quad [ X_3,Z_1 ] =-X_4, \quad [ X_4,Z_1 ] =  X_3, \cr
 & & [ X_1,Z_2 ] =  X_3, \quad [ X_2,Z_2 ] = X_4, \quad [ X_3,Z_2 ] =-X_1, \quad [ X_4,Z_2 ] = -X_2, \cr
 & & [ X_1,Z_3 ] =  X_4, \quad [ X_2,Z_3 ] =-X_3, \quad [ X_3,Z_3 ] = X_2, \quad [ X_4,Z_3 ] = -X_1,
\end{eqnarray}
\begin{eqnarray}
\nonumber
 & & [ Z_1,Z_2 ] = 2Z_3, \quad [ Z_1,Z_3 ] =-2Z_2, \quad [ Z_2,Z_3 ] = 2Z_1.
\end{eqnarray}
The adjoint action of $\fh$ on $\fm$ is given by the operators
\begin{eqnarray}
\label{a33}
{\rm ad}(H_1)|_\fm = A_{12}+A_{34},\quad
{\rm ad}(H_2)|_\fm = A_{13}-A_{24},\quad
{\rm ad}(H_3)|_\fm = A_{14}+A_{23}.
\end{eqnarray}
There is a 7-parameter family of invariant scalar products on $\fm$ and consequently invariant Riemannian
metrics on $M=G/H$.
The geodesic orbit property of these Riemannian metrics was studied
also with respect to extended groups in
\cite{Ni2}, \cite{Z1}, \cite{Z2}, \cite{Z3}, see also \cite{BN} for the overview.
We are now going to generalize this to Finsler $(\alpha,\beta)$ metrics and
formulate the results in terms of geodesic graphs, in the following sections.
Without loosing any Riemannian geodesic orbit metrics, we start with
the $4$-parameter family of invariant scalar products on $\fm$.
Up to a multiple, the basis $B$ will be orthogonal with
$\langle E_i,E_i\rangle = 1$, $\langle Z_1,Z_1\rangle = c_1$,
$\langle Z_2,Z_2\rangle = c_2$, $\langle Z_3,Z_3\rangle = c_3$.
We obtain the 4-parameter family of invariant Riemannian metrics on $G/H={\mathrm{Sp}}(2)/{\mathrm{Sp}}(1)$.
We use again equation (\ref{goleman2}) in the form
\begin{eqnarray}
\label{forex3}
\alpha \Bigl ( X , [\xi(X),U]_{\mathfrak m}\Bigr ) & = &
- \alpha \Bigl ( X , [X,U]_{\mathfrak m}\Bigr ) 
- \alpha \Bigl ( \zeta(X)\cdot  V, [X,U]_{\mathfrak m}\Bigr ) .
\end{eqnarray}

\noindent
We write each vector $X\in\fm$ and each vector $F\in\fh$ in the forms
\begin{eqnarray}
\nonumber
X & = & x_1X_1 + \dots + x_4 X_4 + z_1 Z_1 + \dots + z_3 Z_3, \cr
F & = & \xi_1 H_1 + \dots +\xi_3 H_3
\end{eqnarray}
and consider formula (\ref{forex3}), without the last Finslerian term.
We obtain the system of linear equations whose extended matrix
$(\bA | \bB)$, without the Finslerian part $ (\bC) $, is
{\small
\begin{eqnarray}
\nonumber
\left (
\begin{array}
{ccc|c}
{x_2}&{x_3}&{x_4}&
(1-2c_1) z_1x_2 +(1-2c_2) z_2x_3 +(1-2c_3) z_3x_4
\\{\medskip}
{-x_1}&{-x_4}&{ x_3}&
-(1-2c_1)z_1x_1 +(1-2c_2) z_2x_4 -(1-2c_3) z_3x_3
\\{\medskip}
{ x_4}&{-x_1}&{-x_2}&
-(1-2c_1) z_1x_4 -(1-2c_2) z_2x_1 +(1-2c_3) z_3x_2
\\{\medskip}
{-x_3}&{ x_2}&{-x_1}&
(1-2c_1) z_1x_3 -(1-2c_2) z_2x_2 -(1-2c_3) z_3x_1
\\{\medskip}
0 & 0 & 0 & 2z_2z_3(c_3-c_2) 
\\{\medskip}
0 & 0 & 0 & 2z_1z_3(c_1-c_3)
\\{\medskip}
0 & 0 & 0 & 2z_1z_2(c_2-c_1)
\end{array}
\right ).
\end{eqnarray}
}
We see immediately that the system is solvable (and $G/H$ is a g.o. space)
for the Riemannian metric (with the first right-hand side $\bB$ only) if and only if $c_1=c_2=c_3$.
We see that just the corresponding $2$-parameter family of invariant Riemannian metrics
on ${\mathrm{Sp}}(2)/{\mathrm{Sp}}(1)$ are geodesic orbit metrics.
\begin{prop}
\label{notgo}
For each metric from the $2$-parameter family of geodesic orbit metrics
on the sphere $S^{7}={\mathrm{Sp}}(2)/{\mathrm{Sp}}(1)$,
there exist a three-dimensional ${\mathrm{Ad}}(H)$-invariant subspace in $\fm$.
For each vector $V$ from this subspace, the Riemannian geodesic orbit metric $\alpha$
can be modified using this vector to obtain invariant Finsler $(\alpha,\beta)$ metrics.
None of these purely Finsler metrics is geodesic orbit with respect to the group ${\mathrm{Sp}}(2)$.
\end{prop}
{\it Proof.}
The existence of the $2$-parameter family of Riemannian geodesic orbit metrics is well known,
see for example \cite{BN}, and they are explicitly described in the example above.
From the Lie brackets and from the adjoint representation given by the operators (\ref{a33})
we see that there is the three-dimensional ${\mathrm{Ad}}(H)$-invariant subspace in $\fm$,
generated by the vectors $V_1,\dots,V_3$.
Hence, any nonzero vector $V=v_1 Z_1 + \dots + v_3 Z_3$,
together with a function $\phi(s)$, determines
invariant Finsler $(\alpha,\beta)$ metrics on $M=G/H$, according to Lemma \ref{ablemma}.
Concerning the geodesic orbit condition, we apply the condition $c_1=c_2=c_3$ into the calculations
above and we obtain the extended matrix $(\bA | \bB | \bC)$, also with the Finslerian part $(\bC)$, in the form
{\small
\begin{eqnarray}
\nonumber
\left (
\begin{array}
{ccc|c|c}
{x_2}&{x_3}&{x_4}& 
(1-2c)(z_1x_2 + z_2x_3 + z_3x_4)
& 2\zeta c(- v_1x_2 - v_2x_3 - v_3x_4)
\\{\medskip}
{-x_1}&{-x_4}&{ x_3}&
(1-2c)(-z_1x_1 + z_2x_4 - z_3x_3)
& 2\zeta c( v_1x_1 - v_2x_4 + v_3x_3)
\\{\medskip}
{ x_4}&{-x_1}&{-x_2}&
(1-2c)(-z_1x_4 - z_2x_1 + z_3x_2)
& 2\zeta c( v_1x_4 + v_2x_1 - v_3x_2)
\\{\medskip}
{-x_3}&{ x_2}&{-x_1}&
(1-2c)(z_1x_3 - z_2x_2 - z_3x_1)
& 2\zeta c(-v_1x_3 +v_2x_2 + v_3x_1)
\\{\medskip}
0 & 0 & 0 & 0 &
2 \zeta c ( z_2 v_3 - z_3 v_2 )
\\{\medskip}
0 & 0 & 0 &  0 &
2 \zeta c ( z_3 v_1 - z_1 v_3 )
\\{\medskip}
0 & 0 & 0 &  0 &
2 \zeta c ( z_1 v_2 - z_2 v_1 )
\end{array}
\right ).
\end{eqnarray}
}
We see that for nonzero vector $V\in\fm$ and for the general choice $U=u_1Z_1+\dots+u_3Z_3$,
the system of equations corresponding to the above extended matrix is not solvable for general $X\in\fm$
and hence $(G/H,F)$ is not a geodesic orbit space.
$\hfill\square$

Proposition \ref{notgo} gives the first example to the open question posed in \cite{DuNew}.
However, $(G/H,F)$ becomes geodesic orbit space with respect to the bigger group $\tilde G\supset G$,
according to Theorem \ref{t2}.
We shall continue with the analysis of this situation in more detail.

\subsection{$S^{4n+3}={\mathrm{Sp}}(n+1)\cdot{\mathrm{U}}(1)/{\mathrm{Sp}}(n)\cdot{\mathrm{diag}}({\mathrm{U}}(1))$ }
\label{35}
In the previous section we have seen that any vector of the form $V=v_1Z_1+\dots+v_3Z_3$ is
${\mathrm{Ad}}(H)$-invariant and it determines an invariant Finsler $(\alpha,\beta)$ metric.
For each such metric, the isotropy algebra can be extended by the formal operator ${\mathrm{ad}}(V)$
according to Proposition \ref{prop1} and its proof in \cite{DuNew}
and we obtain $\widetilde\fh\simeq{\mathrm{sp(1)}}\oplus {\mathrm{u(1)}}$.
For the particular example, let us choose $V=v_1Z_1, v_1\in \bR$. Then we have
$\widetilde\fh={\mathrm{span}}\{\fh,W_1\}\simeq{\mathrm{sp(1)}}\oplus {\mathrm{u(1)}}$ for the operator
\begin{eqnarray}
\label{W1}
W_1={\rm ad}(Z_1)|_\fm = 2B_{23}-A_{12}+A_{34}
\end{eqnarray}
and we have the new expression of our manifold in the form
$S^{7}={\mathrm{Sp}}(2)\cdot{\mathrm{U}}(1)/{\mathrm{Sp}}(1)\cdot{\mathrm{diag}}({\mathrm{U}}(1))$.
We will now consider the $3$-parameter family of Riemannian metrics corresponding to
scalar products from previous section, which satisfy the condition $c_2=c_3$ and which are invariant
with respect to $\widetilde\fh$. Related invariant Minkowski norms and Finsler $(\alpha,\beta)$ metrics
are determined by the vector $V$ above and a function $\phi(s)$, according to Lemma \ref{ablemma}.
Formula (\ref{forex3}) gives us the system of linear equations whose extended matrix $(\bA | \bB | \bC)$ is
{\small
\begin{eqnarray}
\nonumber
\left (
\begin{array}
{cccc|c|c}
{x_2}&{x_3}&{x_4}& {-x_2}&
(1-2c_1)z_1x_2 + (1-2c_2)(z_2x_3 + z_3x_4)
& -2\zeta x_2 v_1 c_1
\\{\medskip}
{-x_1}&{-x_4}&{ x_3}& { x_1}&
-(1-2c_1)z_1x_1 + (1-2c_2)(z_2x_4 - z_3x_3)
& 2\zeta x_1 v_1 c_1
\\{\medskip}
{ x_4}&{-x_1}&{-x_2}& { x_4}&
-(1-2c_1)z_1x_4 + (1-2c_2)(- z_2x_1 + z_3x_2)
& 2\zeta x_4 v_1 c_1
\\{\medskip}
{-x_3}&{ x_2}&{-x_1}& {-x_3}&
(1-2c_1)z_1x_3 - (1-2c_2)(z_2x_2 + z_3x_1)
& - 2\zeta x_3 v_1 c_1
\\{\medskip}
0 & 0 & 0 &  2z_3c_2 & 2z_1z_3(c_1-c_2) &
2 \zeta c_1 z_3 v_1  
\\{\medskip}
0 & 0 & 0 & -2z_2c_2 & 2z_1z_2(c_2-c_1) &
- 2 \zeta c_1 z_2 v_1 
\end{array}
\right ).
\end{eqnarray}
}
The rank of this system is equal to four and we see that it is solvable.
Hence the Finsler $(\alpha,\beta)$ metrics determined by $V=v_1Z_1$ are geodesic orbit metrics.
However, for the compoments of the Riemannian geodesic graph $\xi^R$, after cancelling out the common factor,
we have $\xi_i^R=P_i/P$, where ${\mathrm{deg}}(P_i)-1={\mathrm{deg}}(P)=2$
and it is not easy to describe the Finslerian geodesic graph $\xi$ using $\xi^R$.
To do so, we change the reductive decomposition, according to Proposition \ref{prop1}.
We put $\bar Z_1=Z_1-W_1$ (hence $\bar Z_1\in{\mathfrak{c}}(\fg)$ and $Z_1=\bar Z_1+W_1$)
and $\fm'={\mathrm{span}}\{X_i,\bar Z_1,Z_2,Z_3\}$. It holds
\begin {eqnarray}
\nonumber
 [ H_1,H_2 ] =  2 H_3, \quad [ H_1,H_3 ] = -2 H_2, \quad [ H_2,H_3 ] = 2 H_1, \quad [ H_i,W_1 ] = 0
\end{eqnarray}
and the projections of the Lie brackets of generators of the space $\fm'$ onto $\fm'$ are
\begin{eqnarray}
\nonumber
 & & [ X_1,X_2 ]_{\fm'} = -2 \bar Z_1 , \cr
 & & [ X_1,X_3 ]_{\fm'} = -2 Z_2 , \quad [ X_2,X_3 ]_{\fm'} =  2 Z_3 , \cr
 & & [ X_1,X_4 ]_{\fm'} = -2 Z_3 , \quad [ X_2,X_4 ]_{\fm'} =  -2 Z_2 ,
 \quad [ X_3,X_4 ]_{\fm'} = ~2 \bar Z_1 ,
\end{eqnarray}
\begin{eqnarray}
\nonumber
 & & [ X_1,\bar Z_1 ]_{\fm'} =  0,~~ \quad [ X_2,\bar Z_1 ]_{\fm'} =0,~~~ \quad [ X_3,\bar Z_1 ]_{\fm'} =0,~~~
 \quad [ X_4,\bar Z_1 ]_{\fm'} =  0,~~ \cr
 & & [ X_1,Z_2 ]_{\fm'} =  X_3, \quad [ X_2,Z_2 ]_{\fm'} = X_4, \quad [ X_3,Z_2 ]_{\fm'} =-X_1,
 \quad [ X_4,Z_2 ]_{\fm'} = -X_2, \cr
 & & [ X_1,Z_3 ]_{\fm'} =  X_4, \quad [ X_2,Z_3 ]_{\fm'} =-X_3, \quad [ X_3,Z_3 ]_{\fm'} = X_2,
 \quad [ X_4,Z_3 ]_{\fm'} = -X_1,
\end{eqnarray}
\begin{eqnarray}
\nonumber
 & & [ \bar Z_1,Z_2 ]_{\fm'} = 0, \quad [ \bar Z_1,Z_3 ]_{\fm'} = 0, \quad [ Z_2,Z_3 ]_{\fm'} = 2\bar Z_1.
\end{eqnarray}
The adjoint action of $\widetilde\fh$ on $\fm'$ is given by the operators
\begin{eqnarray}
\nonumber
& & {\rm ad}(H_1)|_\fm = A_{12}+A_{34},\quad
{\rm ad}(H_2)|_\fm = A_{13}-A_{24},\cr
& & {\rm ad}(H_3)|_\fm = A_{14}+A_{23}, \quad
{\rm ad}(W_1)|_\fm = 2B_{23}-A_{12}+A_{34}.
\end{eqnarray}
Geodesic lemma and formula (\ref{forex3}) gives us the system of equations
whose extended matrix $(\bA,\bB,\bC)$ is equivalent to
{\small
\begin{eqnarray}
\nonumber
\left (
\begin{array}
{cccc|c|c}
{x_2}&{x_3}&{x_4}& {-x_2}&
-2c_1 z_1x_2 + (1-2c_2)(z_2x_3 + z_3x_4)
& -2\zeta v c_1 x_2 
\\{\medskip}
{-x_1}&{-x_4}&{ x_3}& { x_1}&
2c_1 z_1x_1 + (1-2c_2)(z_2x_4 - z_3x_3)
& 2\zeta v c_1 x_1
\\{\medskip}
{ x_4}&{-x_1}&{-x_2}& { x_4}&
2c_1 z_1x_4 + (1-2c_2)(- z_2x_1 + z_3x_2)
& 2\zeta v c_1 x_4
\\{\medskip}
0 & 0 & 0 &  c_2 & z_1 c_1 &
 \zeta v c_1
\end{array}
\right ).
\end{eqnarray}
}
Using the Cramer's rule, we obtain components of geodesic graph
\begin{eqnarray}
\nonumber
\xi_1 & = & \Bigl [
  (\frac{c_1}{c_2}-2c_1) ( x_1^2 + x_2^2 - x_3^2 - x_4^2 ) (z_1+\zeta v) +\cr
&& 2 (1-2c_2) (x_2x_3 - x_1x_4) z_2 + 2 (1-2c_2) (x_1x_3 + x_2x_4) z_3 \Bigr ]/\|x\|^2 ,\cr
\xi_2 & = & \Bigl [
 2 (\frac{c_1}{c_2}-2c_1) (x_2x_3 + x_1x_4) (z_1+\zeta v) +\cr
&& (1-2c_2) ( x_1^2 - x_2^2 + x_3^2 - x_4^2 ) z_2 + 2 (1-2c_2) (x_3x_4 - x_1x_2) z_3 \Bigr ]/\|x\|^2 ,\cr
\xi_3 & = & \Bigl [
 2 (\frac{c_1}{c_2}-2c_1) (x_2x_4 - x_1x_3) (z_1+\zeta v) +\cr
&& 2 (1-2c_2) (x_1x_2 + x_3x_4) z_2 + (1-2c_2)  ( x_1^2 - x_2^2 - x_3^2 + x_4^2 ) z_3 \Bigr ]/\|x\|^2 ,\cr
\xi_4 & = &
\frac{c_1}{c_2}(z_1 + \zeta v),
\end{eqnarray}
where we put
\begin{eqnarray}
\nonumber
\|x\|^2 = x_1^2 + x_2^2 + x_3^2 + x_4^2.
\end{eqnarray}
This geodesic graph satisfies the statement of Theorem \ref{t2}.
\begin{prop}
\label{pr2}
Each of the $3$-parameter family of geodesic orbit metrics on the sphere
$S^{4n+3}={\mathrm{Sp}}(n+1)\cdot{\mathrm{U}}(1)/{\mathrm{Sp}}(n)\cdot{\mathrm{diag}}({\mathrm{U}}(1))$ 
can be modified in one direction to obtain invariant Finsler $(\alpha,\beta)$ metric.
Each of these Finsler metrics is geodesic orbit metric with respect to the group ${\mathrm{Sp}}(n+1)\cdot{\mathrm{U}}(1)$.
\end{prop}
{\it Proof.}
The existence of the $3$-parameter family of geodesic orbit Riemannian metrics on these spheres is well known,
see for example \cite{BN}.
The detailed construction of these metrics and also the statement of the proposition
in the particular case $n=1$ is clear from the example above.
In any reductive decompositon $\fg=\fh+\fm$, the isotropy representation admits just one one-dimensional
invariant subspace, which is the projection of the center ${\mathfrak{c}}(\fg)$ into $\fm$.
The reductive decomposition can be changed for the decomposition satisfying Proposition \ref{prop1},
which also the example above illustrates.
Geodesic graph for any invariant Finsler $(\alpha,\beta)$ metric $F$ determined by the Riemannian
metric $\alpha$, the one-form $\beta$ which is $\alpha$-equivalent with $V$ and a function $\phi(s)$,
according to Lemma \ref{ablemma}, can be constructed using Theorem \ref{t2}.
This proves the geodesic orbit property of $(G/H,F)$.
$\hfill\square$

\subsection{$S^{4n+3}={\mathrm{Sp}}(n+1)\cdot{\mathrm{Sp}}(1)/{\mathrm{Sp}}(n)\cdot{\mathrm{diag}}({\mathrm{Sp}}(1))$ }
\label{34}
We finish with considering the $2$-parameter family of metrics
determined by the condition $c_1=c_2=c_3$, from the Section \ref{33}.
In such a case, the isotropy algebra can be extended by the operators
\begin{eqnarray}
\label{b33}
&& W_1={\rm ad}(Z_1)|_\fm = 2B_{23}-A_{12}+A_{34},\cr
&& W_2={\rm ad}(Z_2)|_\fm =-2B_{13}-A_{13}-A_{24},\cr
&& W_3={\rm ad}(Z_3)|_\fm = 2B_{12}-A_{14}+A_{23}
\end{eqnarray}
and we obtain $\fh'={\mathrm{span}}(\fh,W_1,\dots,W_3)\simeq{\mathrm{sp(1)}}\oplus {\mathrm{sp(1)}}$
and $\fg'=\fh'+\fm$. We have the new expression of our manifold $M$ in the form
$M=G'/H'={\mathrm{Sp}}(2)\cdot{\mathrm{Sp}}(1)/{\mathrm{Sp}}(1)\cdot{\mathrm{diag}}({\mathrm{Sp}}(1))$.
With respect to the bigger group $G'$, for the vector $F\in\fh'$, we write
\begin{eqnarray}
\nonumber
F & = & \xi_1 H_1 + \dots +\xi_3 H_3 + \xi_4 W_1 + \xi_5 W_1 + \xi_6 W_3
\end{eqnarray}
and the Riemannian version of formula (\ref{forex3}), without the last Finslerian term, gives us the system of equations
whose extended matrix $(\bA | \bB )$, without the Finslerian part $(\bC)$, is
{\small
\begin{eqnarray}
\nonumber
\left (
\begin{array}
{cccccc|c}
{x_2}&{x_3}&{x_4}& {-x_2}&{-x_3}&{-x_4}&
(1-2c)(z_1x_2 + z_2x_3 + z_3x_4)
\\{\medskip}
{-x_1}&{-x_4}&{ x_3}& { x_1}&{-x_4}&{ x_3}&
(1-2c)(-z_1x_1 + z_2x_4 - z_3x_3)
\\{\medskip}
{ x_4}&{-x_1}&{-x_2}& { x_4}&{ x_1}&{-x_2}&
(1-2c)(-z_1x_4 - z_2x_1 + z_3x_2)
\\{\medskip}
{-x_3}&{ x_2}&{-x_1}& {-x_3}&{ x_2}&{ x_1}&
(1-2c)(z_1x_3 - z_2x_2 - z_3x_1)
\\{\medskip}
0 & 0 & 0 & 0 & -2z_3c &  2z_2c & 0
\\{\medskip}
0 & 0 & 0 &  2z_3c & 0 & -2z_1c & 0
\\{\medskip}
0 & 0 & 0 & -2z_2c &  2z_1c & 0 & 0
\end{array}
\right ).
\end{eqnarray}
}
We easily see the linear geodesic graph.
\begin{eqnarray}
\nonumber
\xi(X) & = & (1-2c) \bigl (  -z_1 W_1 - z_2 W_2 - z_3 W_3 \bigr )
\end{eqnarray}
and hence $M=G/H$ is naturally reductive, which is well known, see for example \cite{Z2}, \cite{Z3}.
\begin{prop}
\label{prop12}
None of the $2$-parameter family of naturally reductive metrics on the sphere
$S^{4n+3}={\mathrm{Sp}}(n+1)\cdot{\mathrm{Sp}}(1)/{\mathrm{Sp}}(n)\cdot{\mathrm{diag}}({\mathrm{Sp}}(1))$
can be modified into invariant purely Finsler $(\alpha,\beta)$ metric.
\end{prop}
{\it Proof.}
In the special case $n=1$, we see from the adjoint representation given by the operators
(\ref{a33}) and (\ref{b33}) that there are no ${\mathrm{Ad}}(H')$-invariant vectors in $\fm$
and hence no ${\mathrm{Sp}}(n+1)\cdot{\mathrm{Sp}}(1)$-invariant purely Finsler $(\alpha,\beta)$ metrics.
Concerning the general case, the Lie algebra $\fg'={\mathfrak{sp}}(n+1)\oplus {\mathfrak{sp}}(1)$ is centerless
and the manifold $M=(G'/H',\alpha)$ does not admit expression in an extended form
$M=\widetilde G/\widetilde H$ for $\widetilde G\supset G'$.
According to Corollary \ref{cor1}, there are no ${\mathrm{Ad}}(H')$-invariant vectors in $\fm$
and hence no ${\mathrm{Sp}}(n+1)\cdot{\mathrm{Sp}}(1)$-invariant purely Finsler $(\alpha,\beta)$ metrics.
$\hfill\square$

\subsection{Geodesic orbit $(\alpha,\beta)$ metrics on spheres}
We conclude our observations with the classification theorem.
\begin{theorem}
The only Riemannian geodesic orbit metrics which can be modified into geodesic orbit
purely Finsler $(\alpha,\beta)$ metrics on spheres are:\\
- The two parameter family of metrics on
$S^{2n+1}={\mathrm{U}}(n+1)/{\mathrm{U}}(n)$;\\
- The three parameter family of metrics on
$S^{4n+3}={\mathrm{Sp}}(n+1)\cdot{\mathrm{U}}(1)/{\mathrm{Sp}}(n)\cdot{\mathrm{diag}}({\mathrm{U}}(1))$.\\
In both cases, these modifications can be done in one direction. In a fixed reductive decomposition,
this direction is determined by the projection of the center ${\mathfrak{c}}(\fg)$ into $\fm$.
\end{theorem}
{\it Proof.}
The two parameter family of Riemannian geodesic orbit metrics on
$S^{2n+1}={\mathrm{SU}}(n+1)/{\mathrm{SU}}(n)$
coincides with the two parameter family of Riemannian geodesic orbit metrics on
$S^{2n+1}={\mathrm{U}}(n+1)/{\mathrm{U}}(n)$
and geodesic orbit Finsler $(\alpha,\beta)$ metrics on the later space were described in Proposition \ref{pr1}.
Geodesic orbit Finsler $(\alpha,\beta)$ metrics on the space
$S^{4n+3}={\mathrm{Sp}}(n+1)\cdot{\mathrm{Sp}}(1)/{\mathrm{Sp}}(n)\cdot{\mathrm{diag}}({\mathrm{Sp}}(1))$
were described in Proposition \ref{pr2}.
Riemannian geodesic orbit metrics on spheres
$S^{4n+3}={\mathrm{Sp}}(n+1)\cdot{\mathrm{Sp}}(1)/{\mathrm{Sp}}(n)\cdot{\mathrm{diag}}({\mathrm{Sp}}(1))$
do not admit modification into invariant purely Finsler $(\alpha,\beta)$ metrics,
as we have shown in Proposition \ref{prop12}.
Concerning the exceptional spheres
$S^{6}={\mathrm{G}}_2/{\mathrm{SU}}(3)$,
$S^{7}={\mathrm{Spin}}(7)/{\mathrm{G}}_2$
and $S^{15}={\mathrm{Spin}}(9)/{\mathrm{Spin}}(7)$,
the Lie algebras $\fg_2$, ${\mathfrak{spin}}(7)$ and ${\mathfrak{spin}}(9)$ are centerless.
These spheres with Riemannian geodesic orbit metrics also do not admit presentations in an extended form
$\widetilde G/\widetilde H$ for $\widetilde G\supset G$. Hence, according to Corollary \ref{cor1},
they do not admit ${\mathrm{Ad}}(H)$-invariant vectors in $\fm$ and hence no modification
into invariant purely Finsler $(\alpha,\beta)$ metrics.
$~\hfill\square$

\section{Projective spaces}
As we have done explicit calculations with the spheres in previous section,
it is easy to illustrate similar behaviour with some projective spaces now.
With Riemannian metrics, the geodesic orbit property on these projective spaces was studied
in \cite{Z1}, \cite{Z2}, \cite{Z3}, see also \cite{BN} for the overview.
Again, we shall study explicitly the low-dimensional examples.

\subsection{${\mathbb{C}}P^{n}={\mathrm{SU}}(n+1)/{\mathrm{S}}(\mathrm{U}(n)\cdot{\mathrm{U}}(1))$}
\label{41}
For $n=2$, we have ${\mathbb{C}}P^2={\mathrm{SU}}(3)/{\mathrm{S}}(\mathrm{U}(2)\cdot {\mathrm{U}}(1))$.
We choose the basis $\{H_1,H_2,H_3,Z\}$ of the Lie algebra $\fh={\mathfrak{s}}(\mathfrak{u}(2)\oplus{\mathrm{u}}(1))$,
given by the matrices
\begin{eqnarray}
\nonumber
H_1= \left ( \begin{array} {ccc}
i & 0 & 0
\\{\medskip}
0 &-i & 0
\\{\medskip}
0 & 0 & 0
\end{array} \right ),
H_2= \left ( \begin{array} {ccc}
0 & i & 0
\\{\medskip}
i & 0 & 0
\\{\medskip}
0 & 0 & 0
\end{array} \right ),
H_3= \left ( \begin{array} {ccc}
0 & 1 & 0
\\{\medskip}
-1 & 0 & 0
\\{\medskip}
0 & 0 & 0
\end{array} \right ),
Z= \left ( \begin{array} {ccc}
-\frac{i}{2} & 0 & 0
\\{\medskip}
0 & -\frac{i}{2} & 0
\\{\medskip}
0 & 0 & i
\end{array} \right ).
\end{eqnarray}
For $\fm$, we choose the basis $B=\{X_1,X_2,Y_1,Y_2\}$, given by the matrices in Section \ref{31}.
We use formulas (\ref{z31}) and (\ref{a31}) and we obtain the adjoint action of $\fh$ on $\fm$,
which is given by the operators
\begin{eqnarray}
\nonumber
& {\rm ad}(H_1)|_\fm =  A_{12}-A_{34},\qquad
{\rm ad}(H_2)|_\fm =  A_{14}-A_{23},\cr
& {\rm ad}(H_3)|_\fm = -A_{13}-A_{24},\qquad
{\rm ad}(Z)|_\fm =  -\frac{3}{2} A_{12} -\frac{3}{2}A_{34}.
\end{eqnarray}
We see that there is a $1$-parameter family of invariant Riemannian metrics on $M=G/H$.
These are known to be normal homogeneous and hence naturally reductive and also geodesic orbit spaces.
We also see that there is not any ${\mathrm{Ad}}(H)$-invariant vector $V\in\fm$ and hence
no invariant purely Finsler $(\alpha,\beta)$ metric on $M=G/H$.

\subsection{${\mathbb{C}}P^{2n+1}={\mathrm{Sp}}(n+1)/{\mathrm{Sp}}(n)\cdot {\mathrm{U}}(1)$}
\label{42}
For $n=1$, we have ${\mathbb{C}}P^3={\mathrm{Sp}}(2)/{\mathrm{Sp}}(1)\cdot {\mathrm{U}}(1)$.
We use the matrices from Section \ref{33}. 
Now, we consider the basis $\{H_1,H_2,H_3,Z_1\}$ of $\fh$ and the basis $B=\{X_1,\dots,X_4,Z_1,Z_3\}$ of $\fm$.
Using formulas (\ref{a33}) and (\ref{b33}),
we see that the adjoint action of $\fh$ on $\fm$ is given by the operators
\begin{eqnarray}
\nonumber
& & {\rm ad}(H_1)|_\fm =  A_{12}+A_{34},\qquad
{\rm ad}(H_2)|_\fm =  A_{13}-A_{24},\cr
& & {\rm ad}(H_3)|_\fm =  A_{14}+A_{23},\qquad
{\rm ad}(Z_1)|_\fm = 2B_{23}-A_{12}+A_{34}.
\end{eqnarray}
We see that there is a $2$-parameter family of invariant Riemannian metrics on $M=G/H$.
These metrics are known to be weakly symmetric, some of them are naturally reductive, see \cite{Z2}, \cite{Z3} for details.
The detailed analysis of these Riemannian metrics and geodesic vectors can be found in \cite{BNN}.
We also see that there is not any ${\mathrm{Ad}}(H)$-invariant vector $V\in\fm$ and hence
no invariant purely Finsler $(\alpha,\beta)$ metric on $M=G/H$.

\subsection{${\mathbb{H}}P^n={\mathrm{Sp}}(n+1)/{\mathrm{Sp}}(n)\cdot {\mathrm{Sp}}(1)$}
\label{43}
For $n=1$, we have ${\mathbb{H}}P^1={\mathrm{Sp}}(2)/{\mathrm{Sp}}(1)\cdot {\mathrm{Sp}}(1)$.
We use again the matrices from Section \ref{33}.
Now, we consider the basis $\{H_1,H_2,H_3,Z_1,\dots,Z_3\}$ of $\fh$ and the basis $B=\{X_1,\dots,X_4\}$ of $\fm$.
Using formulas (\ref{a33}) and (\ref{b33}), we see that the adjoint action of $\fh$ on $\fm$ is given by the operators
\begin{eqnarray}
\nonumber
& {\rm ad}(H_1)|_\fm = A_{12}+A_{34},\quad
{\rm ad}(H_2)|_\fm = A_{13}-A_{24},\quad
{\rm ad}(H_3)|_\fm = A_{14}+A_{23}, \cr
& {\rm ad}(Z_1)|_\fm = -A_{12}+A_{34},\quad
{\rm ad}(Z_2)|_\fm = -A_{13}-A_{24},\quad
{\rm ad}(Z_3)|_\fm = -A_{14}+A_{23}.
\end{eqnarray}
We see again that there is a $1$-parameter family of invariant Riemannian metrics on $M=G/H$.
These are known to be normal homogeneous and hence naturally reductive and geodesic orbit spaces.
We also see that there is not any ${\mathrm{Ad}}(H)$-invariant vector $V\in\fm$ and hence
no invariant purely Finsler $(\alpha,\beta)$ metric on $M=G/H$.

\subsection{Geodesic orbit $(\alpha,\beta)$ metrics on projective spaces}
All geodesic orbit metrics on spheres induce geodesic orbit metrics on real projective spaces
${\mathbb{R}}P^n$ in the natural way.
According to \cite{BN}, Riemannian metrics obtained this way and certain Riemannian metrics on
${\mathbb{C}}P^{n}={\mathrm{SU}}(n+1)/{\mathrm{S}}(\mathrm{U}(n)\cdot{\mathrm{U}}(1))$,
${\mathbb{C}}P^{2n+1}={\mathrm{Sp}}(n+1)/{\mathrm{Sp}}(n)\cdot {\mathrm{U}}(1)$,
${\mathbb{H}}P^n={\mathrm{Sp}}(n+1)/{\mathrm{Sp}}(n)\cdot {\mathrm{Sp}}(1)$ and
${\mathbb{C}}aP^2={\mathrm{F}}_4/{\mathrm{Spin}}(9)$
exhaust all Riemannian geodesic orbit metrics on projective spaces.
\begin{prop}
Projective spaces
${\mathbb{C}}P^{n}$, ${\mathbb{C}}P^{2n+1}$, ${\mathbb{H}}P^n$ and ${\mathbb{C}}aP^2$
mentioned above do not admit invariant purley Finsler $(\alpha,\beta)$ metrics.
\end{prop}
{\it Proof.}
For the particular case $n=1$, the statement for the first three type of spaces follows from previous examples.
Concerning the general case, the Lie algebras $\fg$ of the isometry groups $G$ of mentioned spaces are centerless.
These projective spaces with Riemannian geodesic orbit metrics also do not admit presentations in an extended form
$\widetilde G/\widetilde H$ for $\widetilde G\supset G$. Hence, again according to Corollary \ref{cor1},
they do not admit ${\mathrm{Ad}}(H)$-invariant vectors in $\fm$ and hence no modification
into invariant purely Finsler $(\alpha,\beta)$ metrics.
$~\hfill\square$

\section*{Acknowledgements}
The research is supported by grant PID2019-10519GA-C22 funded by $AEI/10.13$ $039/501100011033$.
The first author is also partially supported by grant GR21055 funded by Junta de Extremadura
and Fondo Europeo de Desarrollo Regional.

\end{document}